\documentclass{article}
\usepackage[english]{babel}
\usepackage{amsmath}
\usepackage{amsthm}
\usepackage{amssymb}
\usepackage{amscd}
\usepackage[latin1]{inputenc}

\newtheorem{teo}{Theorem}[section]
\newtheorem{defi}[teo]{Definition}
\newtheorem{lema}[teo]{Lemma}
\newtheorem{prop}[teo]{Proposition}
\newtheorem{coro}[teo]{Corollary}
\newtheorem{obser}[teo]{Remark}
\newtheorem{ejem}[teo]{Example}

\newtheorem{rem}[teo]{Remark}

\def\Hom{\mathop{\rm Hom}\nolimits}

\def\dim {\mathop{\rm dim_k}\nolimits}

\def\Ext {\mathop{\rm Ext}\nolimits}
\def\Tor {\mathop{\rm Tor}\nolimits}

\def\Ker {\mathop{\rm Ker}\nolimits}
\def\Coker {\mathop{\rm Coker}\nolimits}
\def\Im {\mathop{\rm Im}\nolimits}

\def\Alt {\mathop{\rm Alt}\nolimits}
\def\car {\mathop{\rm char}\nolimits}

\begin{document}

\sf
\title{The first cohomology group of the trivial extension of a monomial algebra }
\date{ }

\author{Claude Cibils, Mar\'\i a Julia Redondo and Manuel Saor\'\i n
\thanks{The second author is a research member of CONICET (Argentina).
The third author thanks the D.G.I. of Spain and the Fundaci\'{o}n
"S\'{e}neca" of Murcia for their financial support.}}

\maketitle

\begin{center}
{{\em{ Dedicated to Raymundo Bautista and Roberto
Mart\'\i nez--Villa \\ for their 60th birthday}}}
\end{center}

\begin{abstract}
Given a finite--dimensional monomial algebra $A$ we consider the trivial extension $TA$ and provide formulae,
depending on the characteristic of the field, for the dimensions of the summands $HH_1(A)$ and $\Alt(DA)$ of the
first Hochschild cohomology group $HH^1(TA)$. From these a formula for the dimension of $HH^1(TA)$ can be derived.

\end{abstract}



\section {\sf Introduction }
The purpose of this paper is to study the first Hochschild cohomology group $HH^1(TA)$ of the trivial extension of
a finite--dimensional monomial algebra $A=kQ/<Z>$ where $k$ is a field, $Q$ a finite quiver and $Z$ a set of paths of length at least
two.

Given an algebra $\Lambda$ and a $\Lambda-\Lambda$--bimodule $X$, the Hochschild cohomology groups $H^i(\Lambda,
X)$, introduced in \cite{ho}, are the groups $\Ext^i_{\Lambda -\Lambda}(\Lambda, X)$. In particular if $X=\Lambda$
we write $HH^i(\Lambda)=H^i(\Lambda, \Lambda)$. Analogously, the Hochschild homology groups $H_i(\Lambda ,X)$ are
the groups $\Tor_i^{\Lambda -\Lambda}(\Lambda ,X)$ and we write $HH_i(\Lambda )=H_i(\Lambda ,\Lambda )$. Although
these groups are not easy to compute in general, some approaches have been successful when the algebra $\Lambda$ is
given by a quiver with relations. For instance,  explicit formulae for the dimensions of $HH^i(\Lambda)$ in terms of
those combinatorial data have been found in  \cite{blm,c1,c2,cs,ha}.

The first Hochschild cohomology group plays an important role in the representation theory of algebras since it is
related to the separation properties of the vertices of the quiver of $\Lambda$, and to the notion of (strong)
simple connectedness (see \cite{amjap, ajap, ha2, japs, sk}).  The importance of simply connected algebras follows
from the fact that  we may often reduce the study of indecomposable modules over an algebra to that of the
corresponding simply connected algebras.

Given an algebra $A$ we consider $DA$, the dual $A-A$--bimodule of $A$.  The trivial extension $TA$ is the algebra
whose underlying vector space is $A \oplus DA$, and the product is given by $(a,f)(b,g)=(ab, ag+fb)$ for any $a, b
\in A$, $f, g \in DA$. In other words, $A$ is a subalgebra of $TA$ and $DA$ is a two--sided ideal endowed with the
zero multiplicative structure.

In \cite{cmrs} it has been shown that if $A$ is a finite--dimensional algebra,  $HH^1(TA)$ is a direct sum of four vector spaces.  More precisely,
$$HH^1(TA) = Z(A) \oplus  HH^1(A) \oplus HH_1(A)^*  \oplus \Alt (DA),$$
where $Z(A)$ is the center of $A$ and $\Alt(DA)=\{\varphi\in
\Hom_{A-A}(DA,A): \ \varphi +\varphi^*=0\}$.  As an immediate
consequence we have that $HH^1(TA)$ never vanishes. The first and
second summands have been studied in \cite{cs}, when $A$ is a monomial
algebra, and explicit formulae for their dimensions are
obtained. The Lie algebra structure of $HH^1(TA)$ with respect to
this decomposition is described in \cite{st}.

In this paper we compute the third and the fourth summands for a given monomial algebra $A$ in terms of the
combinatorics of the quiver and of the set of path relations. In particular we provide precise criteria for the
vanishing of those summands.

An important tool for the computations is given by the circuits of a quiver, which are equivalence classes of
cycles under rotation, see Definition \ref{circuit}. Each circuit has a well defined multiplicity; in
characteristic $p$, the circuits which are relevant  have a multiplicity that is not divisible by $p$. We call them $p'$-circuits by
analogy with the $p'$-conjugacy classes of a group.

The paper is organized as follows.  In Section 2 we introduce some notations and definitions, and we provide
formulae for the dimension of the first Hochschild homology group $HH_1(A)$ when the field $k$ has
characteristic zero.  In Section 3 we extend the previous results for any field $k$ of positive characteristic.  In
Section 4 we obtain a dimension formula for the vector space $\Alt(DA)$ and, as an application of the formulae
obtained in the paper, we describe the monomial algebras $A$ for which $HH^1(TA)$ is minimal, that is, $\dim
HH^1(TA)=1$. It turns out that this condition is equivalent to $HH^1(A)=0$. We do not know if this equivalence holds for a
wider class of algebras -- this question has been pointed out by the referee.

All the algebras considered here are finite--dimensional, but results in Section 2 and 3 hold also for infinite--dimensional monomial algebras.

\section {\sf Degree one Hochschild homology of a monomial algebra}\label{homologia}

Let $k$ be a field, $Q$ a finite quiver with set of vertices
$Q_0$, set of arrows $Q_1$ and $s, t : Q_1 \to Q_0$ be the maps
providing each arrow $a$ with its source vertex $s(a)$ and its
terminal vertex $t(a)$.  A path $\alpha$ of length $l$ is a sequence
of $l$ arrows $a_l \dots a_1$ such that $t(a_i)=s(a_{i+1})$.  We
put $s(\alpha)=s(a_1)$ and $t(\alpha)=t(a_l)$.  Any vertex $u$ is
a trivial path of length zero and we put $s(u)=t(u)=u$.  A cycle
is a path $\alpha$ such that $s(\alpha)=t(\alpha)$; vertices are
always cycles. The corresponding path algebra $kQ$ is the vector
space with basis all the paths in $Q$ and whose product on the basis
elements is defined by the concatenation of the sequences of arrows
of the paths $\beta$ and $\alpha$ if they form a path (namely, if
$t(\alpha)=s(\beta)$) and zero otherwise.  Vertices form a
complete set of orthogonal idempotents.  Note that $\alpha
s(\alpha)=\alpha$ while $\alpha u =0$ if $u \not = s(\alpha)$.  We
have also $t(\alpha) \alpha=\alpha$ and  $u \alpha =0$ if $u \not
= t(\alpha)$. Since all summands in $HH^1(TA)$ are additive with
respect to the decomposition of $A$ into a finite direct product
of algebras, it is not restrictive for our purposes to assume that
$Q$ is connected, something that we do from now on.

Next consider $Z$ a set of paths of length at least two which is
minimal with respect to the subpath order relation, namely, for
each $\gamma \in Z,$ strict subpaths of $\gamma$ are not in $Z$. We
denote $<Z>$ the two-sided ideal generated by $Z$, and $kQ/<Z>$ is
by definition a monomial algebra. We fix $Q$ and $Z$ and put
$A=kQ/<Z>$ in the sequel.

The purpose of this section is to provide a combinatorial formula
computing the dimension of the first Hochschild homology vector
space of a monomial algebra, when $\car k$ = $0$.

We need some  notation in order to describe a chain complex computing Hochschild homology.   Let $B$ be the set
of paths of $Q$ which do not contain any path of $Z$.  In other words, $B$ is a basis of a subvector space of $kQ$
complementing $<Z>$ and $B$ is identified with a $k$--basis of the monomial algebra.  Product of paths in $B$ is
given by usual concatenation, which can be zero if it contains a path from $Z$.  Paths containing a path from $Z$
are called zero paths.

Now we describe  the set of {\it cyclic pairs} of paths.  Let $X$ and $Y$ be sets of paths.  Then
$$X \odot Y = \{ (\alpha, \beta) : t(\beta)=s(\alpha) \ \mbox{and}\ t(\alpha)=s(\beta) \}.$$
Note that $Q_0 \odot X$ is just the set of cycles in $X$.

We denote by $k(X \odot Y)$ the vector space with basis the
set $X \odot Y$. This notation is very convenient to describe a
chain complex whose homology is $HH_*(A)$. We denote by $Z_n$
instead of $\Gamma_n$ the set of n-chains defined in \cite{AG}
(hence, $Z_2=Z$, $Z_1=Q_1$ and $Z_0=Q_0$). This set is independent
of whether we consider left or right modules (cf \cite{b}[Lemma 3.1]) and, as proved in
\cite{b}[Theorem 4.1], the minimal projective resolution $P_*$ of
$A$ as an $A-A$--bimodule is defined as follows: $P_n=\oplus_{p^n\in
Z_n}Ae_{t(p^n)}\otimes e_{s(p^n)}A$,  and the differential
$d:P_{n+1}\longrightarrow P_n$ maps $e_{t(p^{n+1})}\otimes
e_{s(p^{n+1})}$ onto $\sum_{1\leq i\leq
r}(-1)^{e_i}\theta_ie_{t(p_i^n)}\otimes e_{s(p_i^n)}\mu_i$,
provided $\{p_1^n, \dots ,p_r^n\}$ is the set of subpaths of $p^{n+1}$ which
are n-chains and $p^{n+1}=\theta_ip_i^n\mu_i$ is the corresponding
(unique) factorization (here $e_i\in\{0,1\}$ is a convenient
exponent introduced in order to avoid distinction between the even
and odd cases of \cite{b}). Now $HH_*(A)$ is the homology of the
chain complex $P_*\otimes_{A-A}A$ and we leave it as an exercise for
the reader to check that $P_n\otimes_{A-A}A$ is isomorphic to
$k(Z_n\odot B)$ as a $k$-vector space. Viewing that isomorphism as
an identification, the following is now straightforward.

\begin{lema}
\label{Bardzell} The Hochschild homology of $A$ is the homology of
the complex having $k(Z_n\odot B)$ in the n-th position and
differential  $d:k(Z_{n+1}\odot B)\longrightarrow k(Z_n\odot B)$
given by $d(p^{n+1},b)=\sum_{1\leq i\leq
r}(-1)^{e_i}(p_i^n,\mu_ib\theta_i)$, if we make the convention
that a summand is zero when $\mu_ib\theta_i$ is a zero path.
\end{lema}

This complex will be referred to as Bardzell's complex
and its initial part is:

$$ \dots \to k(Z\odot B) \stackrel{d_1}{\rightarrow}  k(Q_1 \odot B) \stackrel{d_0}{\rightarrow} k(Q_0 \odot B) \to 0$$
where
\begin{eqnarray*}
d_0(a, \beta) & = & (t(a), a\beta) - (s(a), \beta a) \\
d_1(\alpha, \beta) &=& \sum_{i=1}^n (a_i, a_{i-1} \dots a_1 \beta
a_n \dots a_{i+1})
\end{eqnarray*}
for $\alpha = a_n \dots a_1$ any path in $Z$.

Next we show that the above complex decomposes  along the circuits
of $Q$ that we define below.

\begin{defi}
\label{circuit} Among the set of cycles of a quiver, consider the equivalence relation generated by
$$\gamma_n \dots \gamma_1 \sim \gamma_1 \gamma_n \dots \gamma_2.$$
The second path is called the rotated of the first path.  An equivalence class for this
relation is by definition a circuit, and we denote by $\cal C$ the set of circuits. A circuit is said to be trivial if it corresponds to a
vertex.
\end{defi}

Cyclic pairs of paths provide circuits by concatenation, namely we have a map
\begin{eqnarray*}
X\odot Y &\rightarrow & \cal C \\
(\alpha, \beta) & \mapsto & \overline{\alpha \beta} = \overline {\beta \alpha}.
\end{eqnarray*}
If $C$ is a fixed circuit we denote by $(X\odot Y)_C$ the fiber over $C$ of this map, namely
$$(X\odot Y)_C = \{ (\alpha, \beta) \in X\odot Y : \alpha \beta \in C \}.$$

\begin{prop}
There exists a decomposition $k(Z_n\odot B)=\oplus_{C \in \cal
C} k(Z_n\odot B)_C$ which is preserved by the differentials of
Bardzell's complex. In particular,
$$HH_n(A)= \oplus_{C \in \cal C} HH_{n,C}(A),$$
where $HH_{*,C}(A)$ is the homology of the $C$-split part of the
complex.
\end{prop}

\begin{proof}
The decomposition is clear. If $(p^{n+1},b)\in (Z_{n+1}\odot
B)_C$, then $d(p^{n+1},b)=\sum_{1\leq i\leq
r}(-1)^{e_i}(p_i^n,\mu_ib\theta_i)$ has the property that
$(p_i^n,\mu_ib\theta_i)\in (Z_n\odot B)_C$, due to the fact that
$p^{n+1}=\theta_ip_i^n\mu_i$.

\end{proof}

\begin{obser}
\sf The above decomposition also holds for cyclic homology. In \cite{iz}, relations between cyclic homology and
global dimension for monomial algebras are obtained; see also \cite{ig}.
\end{obser}

We will compute for each circuit $C$  the corresponding first
homology group $H_{1,C}:=HH_{1,C}(A)$.  Of course if $(Q_1\odot
B)_C = \emptyset$ then $H_{1,C}=0$. Hence we concentrate on the
set of circuits $C$ such that $(Q_1\odot B)_C \not = \emptyset$.
In the sequel we will need the following definitions.

\begin{defi}\label{def circuitos}
\begin{itemize}
\item [ ] \item[i)] A circuit $C$ is said to be useful if
$(Q_1\odot B)_C \not = \emptyset$. \item [ii)] A circuit $C$ is
said to be strong if each cycle of $C$ is a basis vector.  In
other words there is no zero cycle belonging to the circuit $C$.
\item [iii)] A circuit $C$ is said to be efficient if it is a
useful (non strong) circuit satisfying $(Z \odot B)_C \not = \emptyset$.
\end{itemize}
\end{defi}

While computing $H_1$ we will need an evaluation of the differences $\vert Q_1\odot B \vert - \vert Q_0 \odot B\vert$.

\begin{defi}\label{doblev}
Let $W$ be the set of cycles in $Q$ containing precisely one path of $Z$ located at its end.  In other words, a
cycle $\gamma$ is in $W$ if $\gamma = \xi \alpha$ with $\xi \in Z$ and no other subpath of $\gamma$ belongs to
$Z$.  Let $w = \vert W \vert$ and $w_C = \vert W \cap C \vert$ for $C$ a given circuit.
\end{defi}

\begin{lema} \label{difference}

\begin{eqnarray*}
\vert Q_1\odot B \vert - \vert Q_0 \odot B\vert & = & w - \vert Q_0 \vert, \  \mbox{and}\\
\vert (Q_1\odot B)_C \vert - \vert (Q_0 \odot B)_C \vert & = & w_C \ \ \mbox{if $C$ is a non trivial circuit.}
\end{eqnarray*}

\end{lema}

\begin{proof}
Let $B^+$ denote  the set $B \setminus Q_0$ of non zero paths of positive length, and consider the map
$$Q_0 \odot B^+ \stackrel{\varphi}{\rightarrow} Q_1 \odot B$$
which removes the last arrow of a path of $B^+$ and inserts it as a first component:
$$\varphi(s, a_n \dots a_1) = (a_n, a_{n-1} \dots a_1).$$
This map is clearly injective.  The complement of its image consists of cyclic pairs $(a, \beta) \in Q_1 \odot
B$ such that $a \beta$ is no longer in $B$.  Since $\beta \in B$, we have that $a \beta$ contains paths from
$Z$ whose last arrow is $a.$
 Moreover, since $Z$ is minimal, $a \beta$ contains exactly one path of
$Z$ located at its end, which means that $a \beta \in W$.

Conversely, each cycle $\gamma \in W$ has positive length and if $\gamma = a_n \dots a_1$ then $(a_n, a_{n-1}
\dots a_1) \in (Q_1 \odot B) \setminus \Im \varphi$. We have proved that $\vert Q_1\odot B \vert - \vert Q_0 \odot
B^+\vert  = w$.  Note that $Q_0 \odot B = (Q_0 \odot B^+) \cup Q_0$ and this provides the complete formula.  The
specialized formula for a non trivial circuit is clear.
\end{proof}

Recall from 2.5 the definition of strong and useful circuits,  and
note that any non trivial strong circuit is useful.

\begin{prop}
Let $C$ be a useful strong circuit, that is, a non trivial strong circuit.  Then $\dim H_{1,C} =1$.
\end{prop}

\begin{proof}
Since $C$ is strong, $(Z \odot B)_C = \emptyset$, otherwise there would be a zero cycle in the strong circuit $C$.
Therefore $\dim H_{1,C} = \dim \Ker d_{0,C}$. Moreover, $C$ being strong also implies $W \cap C = \emptyset$,
consequently the preceding lemma shows that $\vert (Q_0 \odot B)_C \vert = \vert (Q_1 \odot B)_C \vert$. Hence
$\dim H_{1,C} = \dim \Coker d_{0,C}$.  Recall that $d_{0,C}(a, \beta)= a \beta - \beta a$.  Since $C$ is strong,
$a \beta$ and $\beta a$ are in $B$ for each $(a, \beta) \in (Q_1 \odot B)_C$.  So each difference of a cycle of
$C$ with its rotated cycle is in $\Im d_{0,C}$.  In fact the image of the basis elements of $(Q_1 \odot B)_C$ are
precisely those.  Now the square matrix of $d_{0,C}$ is

\[ \left(\begin{array}{cccccc}
1 & 0 & 0 & & -1 \\
-1 & 1 & 0 \\
0 & -1 & 1  \\
& & & \ddots \\
 & & & &  1
\end{array} \right) \]
which shows that $\dim \Coker d_{0,C}=1$.
\end{proof}

We turn now to useful circuits which are not strong  in order to complete the computation of $HH_1$ along the circuits.

\begin{lema}
Let $C$ be a useful circuit  which is not strong.  Then $d_{0,C}$ is surjective.
\end{lema}

\begin{proof}
Consider the  equivalence relation on circuits restricted  to $B\cap C$. Since $C$ is not strong there is at least one
zero cycle in  $C$, which implies that each equivalence class in $B
\cap C$ is now totally ordered, with the elementary
step for this ordering  given by rotation.  The first and the last element of each of those totally ordered
classes are attained by $d_{0,C}$, as well as all the successive differences. This shows that each element of the
class is in the image of $d_{0,C}$.
\end{proof}

\begin{prop}
Let $C$ be a useful circuit  which is not strong,  satisfying $(Z \odot B)_C = \emptyset$. Then $H_{1,C}=0$.
\end{prop}

\begin{proof}
Since $d_{1,C} =0$ we have that $H_{1,C}= \Ker d_{0,C}$.  By the preceding result $d_{0,C}$ is surjective, hence
$$\dim H_{1,C} = \vert (Q_1 \odot B)_C\vert - \vert (Q_0 \odot B)_C \vert =w_C$$
where $w_C$ is the number of cycles in $C$ which have only one subpath from $Z$ located at its end.  But if  there
exists such a path, then $(Z \odot B)_C$ would be non empty.  Hence $w_C=0$.
\end{proof}

We focus now on efficient circuits, see Definition \ref{def
circuitos}.

\begin{prop} \label{w-1}
Let $C$ be an efficient circuit.  If $k$ has  characteristic zero we have $\dim H_{1,C}=w_C-1$.
\end{prop}

In order to prove this result we will define a canonical element
$K_C$ in $k(Q_1 \odot B)_C$ which will be also useful in positive
characteristic.  First notice that if $\gamma = a_m \dots a_1$ is
a cycle in the quiver, each $a_i$ arising in the sequence of
arrows has a well defined complement path in the cycle, namely
$$a_i^{co} = a_{i-1} \dots a_1 a_m \dots a_{i+1}.$$
Moreover $(a_i, a_i^{co})$ is a cyclic pair of paths.  Note that $a_i$ and $a_j$ can coincide as arrows, but in
general $a_i^{co} \not = a_j^{co}$.  Note however that $a_i^{co} = a_j^{co}$ can also occur  for $i \not = j$.

\begin{defi}
Let $C$ be a circuit and $\gamma = a_m \dots a_1$ be any cycle belonging to $C$.  Let $K_C = \sum_{i=1}^m (a_i,
a_i^{co}) \in k(Q_1 \odot B)_C$.
\end{defi}

\begin{obser}
\begin{itemize} \item [ ]
\item [a)] Our convention is in force, namely if for some $i$ we
have $a_i^{co} \not \in B$ then the pair $(a_i, a_i^{co})$ is
considered as zero in $k(Q_1 \odot B)$. \item [b)] It is clear
that $K_C$ does not depend on the choice of the cycle in $C$.
\end{itemize}
\end{obser}

\begin{lema}
If $k$ has characteristic zero and $C$ is a useful cycle, then $K_C \not =0$.
\end{lema}

\begin{proof}
Since $(Q_1 \odot B)_C \not = \emptyset$ there exists at least one arrow $a_i$ of a cycle $\gamma$ of $C$ such
that $a_i^{co} \in B$.  Note that all the coefficients appearing in the definition of $K_C$ are equal to one, and
no sum of them can be zero in characteristic zero.
\end{proof}

\begin{prop}
Let $C$ be an efficient circuit. Then $\Im d_{1,C} = k K_C$.
\end{prop}

\begin{obser}
This result does not depend on the characteristic.   Note however
that in positive characteristic $K_C$ can be zero.
\end{obser}

\begin{proof}
Since $C$ is efficient, there exists a cyclic pair $(\alpha, \beta) \in (Z \odot B)_C$ and we use the cycle
$\beta \alpha = b_n \dots b_1 a_m \dots a_1$ in order to construct $K_C$:
$$K_C= (a_1, \beta a_m \dots a_2) + (a_2, a_1 \beta a_m \dots a_3) + \dots + (a_m, a_{m-1} \dots a_1 \beta).$$
Indeed each term of the form $(b_j, b_{j-1} \dots b_1 \alpha b_n \dots b_{j+1})$ is zero since its second
component is not in $B$.  Note that some of the written terms in $K_C$ can also be zero, but this has no incidence
in this proof. By definition of $d_{1,C}$ we obtain $d_{1,C}(\alpha, \beta) = K_C$.
\end{proof}

\begin{proof} of Proposition \ref{w-1}. The proof  is now obvious since $C$ is efficient, useful and not strong, hence $d_{0,C}$ is surjective.  In characteristic zero we have proved that $\dim \Im d_{1,C} =1$, while $\vert (Q_1 \odot B)_C \vert - \vert (Q_0 \odot B)_C \vert = w_C$.
\end{proof}

The results we have obtained show that in characteristic zero the
contributing circuits for $HH_1$ are the non trivial strong
circuits and the efficient ones.  The following statement is
obtained by assembling the previous results.

\begin{coro} \label{cero}
Let $A=kQ/<Z>$ be a monomial algebra,  with $k$ a field of
characteristic zero and $Q$ a finite connected quiver.  The
following assertions are equivalent:

\begin{enumerate}
\item $HH_1(A)=0.$
\item Every non trivial circuit of $Q$ contains a
zero cycle and,  whenever $C$ is a circuit such that $(Q_1\odot
B)_C\neq\emptyset\neq (Z\odot B)_C$, there is exactly one pair
$(\xi ,\beta )\in (Z\odot B)_C$ such that the cycle $\xi\beta$
contains no zero relation apart from $\xi$.
\end{enumerate}

\end{coro}

More generally, the following formula holds as a direct consequence of
the previous discussion.

\begin{teo}
\label{teorema1} Let $A=kQ/<Z>$ be a monomial algebra,  with $k$ a
field of characteristic zero and $Q$ a finite connected quiver. We
have
$$\dim HH_1(A) = s + \sum_{C \in \cal E} w_C -e $$
where $\cal E$ is the set of efficient circuits, $e = \vert \cal E
\vert$, $s$ is the number of non trivial strong
circuits  and $w_C$ is defined in Definition \ref{doblev}.
\end{teo}

A formula avoiding the integers $w_C$ can also be obtained as follows.

\begin{coro}
Let $A=kQ/<Z>$ be a monomial algebra,  with $k$ a field of
characteristic zero and $Q$ a finite connected quiver.  Then
$$\dim HH_1(A) = \vert Q_1 \odot B \vert - \vert Q_0 \odot B \vert + \vert Q_0 \vert -e + s.$$
\end{coro}

\begin{proof}
We compute $\sum_{C \in \cal E} w_C$ in order to replace it in the previous result.  Note that if $C$ is trivial
then $w_C=0$. We assert that also if $C$ is not efficient then $w_C =0$.  Indeed, if $(Z \odot B)_C = \emptyset$
then clearly $w_C =0$, while if $(Q_1 \odot B)_C = \emptyset$ and $C$ is not trivial then $(Q_0 \odot B)_C =
\emptyset$ since there is an injective map $\varphi_C : (Q_0 \odot
B^+)_C \to (Q_1 \odot B)_C$, as in the proof of
Lemma \ref{difference}.  This implies that $w_C = \vert (Q_1 \odot B)_C \vert - \vert (Q_0 \odot B)_C \vert =0$.
Consequently
$$ \sum_{C \in \cal E} w_C = \sum_{ C \in \cal C} w_C =
\vert Q_1 \odot B \vert - \vert Q_0 \odot B \vert + \vert Q_0 \vert .$$
\end{proof}

\begin{ejem}
The algebra with quiver $Q$ whose set of vertices is
$Q_0=\mathbf{Z}_5$, with arrows $a_i:i\longrightarrow i+1$,  for
$i\in Q_0$, subject to the relations $a_4a_3a_2=0=a_3a_2a_1$, satisfies $HH_1(A)=0$, whereas for the canonical circuit $C$ we have
$2=\vert (Z\odot B)_C \vert
>w_C=1$. This shows that in Corollary \ref{cero} the last condition
 cannot be replaced by $\vert (Z\odot B)_C \vert =1$.
\end{ejem}

\section {\sf Degree one Hochschild homology of a monomial algebra in
positive characteristic}

Let $C$ be a circuit and $\gamma= a_n \dots a_1$ be a cycle in $C$.  Among the possible iterated rotated cycles of
$\gamma$ clearly the $n$-th one coincides with $\gamma$.  Let $l$ be the smallest integer such that the $l$-th
rotated of $\gamma$ coincides with $\gamma$.  We have $l\!\mid\!n$ and $l$ is called the {\it period} of the
circuit while $m= \frac n l$ is its {\it multiplicity}.

Recall that $K_C$ is defined using any cycle $\gamma=a_n \dots a_1$ of the circuit $C$,
$$K_C= \sum_{i=1}^n (a_i, a_i^{co}) \in k(Q_1 \odot B)_C.$$
Let $\Delta_C = \sum_{(a, \beta) \in (Q_1 \odot B)_C} (a, \beta)$.
Clearly if $C$ is useful then $\Delta_C \not = 0$
for $k$ a field of any characteristic since $\Delta_C$ is the sum of all the basis vectors, a non empty set.

The following result follows from the above considerations.

\begin{lema}
Let $C$ be a circuit of multiplicity $m$.  Then $K_C=m \Delta_C$.
\end{lema}

\begin{defi}
Let $p$ be a prime number.  A circuit is called a $p'$-circuit if its multiplicity is not divisible by $p$.
\end{defi}

\begin{prop}
Let $C$ be an efficient circuit.  If $C$ is a $p'$-circuit then $\dim H_{1,C}=w_C-1$ and $\dim H_{1,C} = w_C$
otherwise.
\end{prop}

\begin{proof}
For an efficient circuit we have proved that $\Im d_{1,C}=k K_C=km \Delta_C$.  If $C$ is a $p'$-circuit, $\dim \Im
d_{1,C} =1$. Otherwise $d_{1,C}=0$.
\end{proof}

The next result computes the dimension of $HH_1(A)$ for a field
$k$ of positive characteristic using the same decomposition as in
characteristic zero but considering $p'$-circuits.

\begin{teo}
\label{char-p} Let $A=kQ/<Z>$ be a monomial algebra,  with $k$ a
field of characteristic $p>0$ and $Q$ a finite connected quiver,
and let $e_{p'}$ be the number of efficient $p'$-circuits in $Q$.
Then
\begin{eqnarray*}
\dim HH_1(A) & = & s + \sum_{C \in \cal E} w_C - e_{p'} \\
& =& \vert Q_1 \odot B \vert - \vert Q_0 \odot B \vert + \vert Q_0 \vert - e_{p'} + s,
\end{eqnarray*}
where $\cal E$ is the set of efficient circuits, $s$ is the number of non trivial strong circuits and $w_C$ is
defined in Definition \ref{doblev}.

\end{teo}

\begin{coro}
\label{cero-p} Let $A=kQ/<Z>$ be a monomial algebra,  with $k$ a
field of characteristic $p>0$ and $Q$ a finite connected quiver.
Then $HH_1(A)=0$ if and only if every non trivial circuit contains a
zero cycle and every circuit $C$ with $(Q_1\odot
B)_C\neq\emptyset\neq (Z\odot B)_C$ is a $p'$-circuit such that
there is exactly one pair $(\xi ,\beta )\in (Z\odot B)_C$ for
which the cycle $\xi\beta$ contains no zero relation apart from
$\xi$.
\end{coro}

\section {\sf The vector space $\Alt(DA)$}\label{alt}

Recall that $\Alt(DA)=\{ \varphi \in \Hom_{A-A}(DA,A) : \varphi +
\varphi^* =0 \}$, where $\varphi^*:DA\longrightarrow DDA\cong A$
is the transpose of $\varphi$ . We start this section by
describing a basis for the vector space $\Hom_{A-A}(DA,A)$. Since
by adjointness we have
$$\Hom_{A-A}(DA,A)= (DA \otimes_{A-A} DA)^*,$$
we describe first a set of generators of $DA \otimes_{A-A} DA$, which will allow us to find the desired basis.

Recall that $B$ is the set of paths in $Q$ which do not contain
any path of $Z$.  The dual basis $B^*$ is a basis of the vector
space $DA$ whose $A-A$--bimodule  structure is given by
$(afb)(x)=f(bxa)$ for any $a,b,x \in A$, $f \in DA$.  This means
that for any $\alpha, \beta, \gamma \in B$, $\alpha \gamma^* \beta
\not =0$ if and only if $\gamma= \beta \xi \alpha$, for some $\xi
\in B$.  In this case, $\alpha \gamma^* \beta =\xi^*$. In
particular, $u\gamma^*v=(v\gamma u)^*$ for any $u, v \in Q_0$,
$\gamma \in B$. This implies that the set of cyclic pairs of
paths
$$B\odot B= \{ (\alpha, \beta) : t(\beta)=s(\alpha) \ \mbox{and}\ t(\alpha)=s(\beta) \}$$
provides a set of generators for $DA \otimes_{A-A} DA$, that is, the set $\{ \alpha^* \otimes \beta^* :  (\alpha,
\beta) \in B\odot B \}$.

\begin{defi}
A cyclic pair $(\alpha, \beta) \in B\odot B$ is said to be neat if the following conditions hold:
\begin{itemize}
\item [i)] if $\beta a \in B$ ($a \beta \in B$) for some $a \in Q_1$ then $a$ is the last (first) arrow in $\alpha$;
\item [ii)] if $\alpha b \in B$ ($b \alpha \in B$) for some $b \in Q_1$ then $b$ is the last (first) arrow in $\beta$.
\end{itemize}
\end{defi}

We denote by $\sim$ the equivalence relation on $B\odot B$
generated by the elementary relations
\begin{eqnarray*}
(a \alpha, \beta) & \sim & (\alpha, \beta a) \\
(\alpha b, \beta) & \sim & (\alpha, b \beta), \ a, b \in Q_1.
\end{eqnarray*}

An equivalence class for this relation will be called {\it neat}
when all its elements are neat cyclic pairs. Let $\cal N $ be the
set of neat equivalence classes.

\begin{lema} \label{gen}
There is a bijective map between  $\cal N$ and  a set of
generators of the vector space $DA \otimes_{A-A} DA$.
\end{lema}

\begin{proof}
We know that $B \odot B$ provides a set of generators. We assert that $\alpha^* \otimes \beta^* =0$ if $(\alpha,
\beta)$ is not neat, and, $\alpha^* \otimes \beta^* = \gamma^* \otimes \delta^*$ if $(\alpha, \beta) \sim (\gamma,
\delta)$. Clearly these facts will prove the Lemma.

Let $(\alpha, \beta)$ be a cyclic pair which is not neat.  We may assume that there exists an arrow $a \in Q_1$
such that $\beta a \in B$ and $a$ is not the last arrow in $\alpha$.  Then
$$\alpha^* \otimes \beta^* = \alpha^* \otimes a(\beta a)^* = \alpha^* a \otimes (\beta a)^*=0.$$
The proof for the other cases is analogous.

To finish the proof we have to check the asserted equality for the elementary relation used to define $\sim$. If
$(a \alpha, \beta) \sim (\alpha, \beta a)$ then
$$(a \alpha)^* \otimes \beta^* = (a \alpha)^* \otimes a( \beta a)^* = \alpha^* \otimes (\beta a)^*.$$
\end{proof}

We are now in a position to describe a basis for
$\Hom_{A-A}(DA,A)$. Let $N\in\cal N$ be a neat equivalence class
and consider the map $\psi_N \in \Hom_k(DA,A)$ defined by
$$\psi_N(\gamma^*)= \sum_{(\gamma, \delta)\in N} \delta$$
for any $\gamma^* \in B^*$.

\begin{lema}
The map $\psi_N$  is a morphism of $A-A$--bimodules.
\end{lema}

\begin{proof}
We shall prove only that $\psi_N$ is a morphism of left
$A$-modules, since the proof that it is a morphism of  right $A$--modules is similar.

It is clear that $\psi_N(u \gamma^*) = u \psi_N(\gamma^*)$ for any
$u \in Q_0$, $\gamma \in B$. In order to finish the proof we have
to see that $\psi_N(a \gamma^*) = a \psi_N(\gamma^*)$ for any $a
\in Q_1$, $\gamma \in B$. Suppose first that $a \gamma^*=0$. This
means that $a$ is not the first arrow in $\gamma$.  If
$\psi_N(\gamma^*)=0$ we are done.  If not, $a \psi_N(\gamma^*)=
\sum_{(\gamma, \delta)\in N} a \delta = 0$ because the cyclic
pairs $(\gamma, \delta)$ are neat.  Now suppose that $a \gamma^*
\not =0$.  Hence $\gamma= \xi a \in B$ and $a \gamma^*=\xi^*$.  It
is clear that $(\xi ,\mu )\in N$ if and only if $\mu =a\delta$, for some
$\delta\in B$ such that $(\gamma ,\delta )\in N$. Then
$$\psi_N(a \gamma^*) = \psi_N(\xi^*)=
\sum_{(\xi, \mu)\in N} \mu = \sum_{(\xi a, \delta)\in N} a \delta
\ = \ \ a \sum_{(\gamma, \delta) \in N} \delta =
a\psi_N(\gamma^*).$$
\end{proof}

\begin{prop}
\label{basis} The set $\{ \psi_N  \}_{N \in \cal N}$ is a basis of
the vector space $\Hom_{A-A}(DA,A)$.
\end{prop}

\begin{proof}
By the adjunction isomorphism
$$\theta : \Hom_{A-A}(DA,A) \to (DA \otimes_{A-A} DA)^*$$ and Lemma
\ref{gen} we conclude that $\{ \psi_N  \}_{N \in \cal N}$ is a
generating set. Considering now the canonical bases $B^*$ and $B$
of $DA$ and $A$, we can identify every $k$-linear map $\varphi
:DA\longrightarrow A$ with its associated matrix, i.e., with the
map $\tilde{\varphi}:B\times B\longrightarrow k$ determined by $\varphi (\alpha^*)=\sum_{\beta\in B}\tilde{\varphi}(\alpha
,\beta )\beta$. Then  $\tilde{\psi}_N$ is the characteristic
function of $N$, i.e., $\tilde{\psi}_N(\alpha ,\beta )=1$ if
$(\alpha ,\beta )\in N$ and $0$ otherwise. From that the
$k$-linear independence of the $\psi_N$'s follows.

\end{proof}

We are now in a position to compute the dimension of the vector
space $\Alt(DA)$.

\begin{lema}
\label{involution} The set $\cal N$ has an involution provided by
the flip of pairs.
\end{lema}

\begin{proof}
Observe that the flip of a cyclic pair provides a cyclic pair.
Moreover neat pairs are preserved and the flip is compatible
with the equivalence relation, namely it is clear that the flips
of elementary equivalent pairs provide elementary equivalent
pairs.
\end{proof}

\begin{teo}
\label{teorema3} Let $A=kQ/<Z>$ be a monomial algebra, where $Q$
is a finite connected quiver. Then
$$\dim \Alt (DA) =
\begin{cases}
\frac{r-s}{2} \ & \mbox{if $\car k \not =2$} \\
\frac{r+s}{2} \ & \mbox{if $\car k =2$}
\end{cases}$$
where $r$ is the number of neat equivalence classes and $s$ is the
number of symmetric ones.
\end{teo}

\begin{proof}
With the same terminology as in the proof of Proposition
\ref{basis}, if $\sigma$ is the above mentioned involution of
$\cal N$, we clearly have $\tilde{\psi}_{\sigma (N)}(\alpha ,\beta
)=\tilde{\psi}_N (\beta ,\alpha )$  for all $(\alpha ,\beta
)\in\sigma (N)$, and hence $\psi_{\sigma (N)}=\psi_N^*$. Therefore a
generic element $\varphi =\sum_{N\in\cal N}\lambda_N\psi_N$ is in
$\Alt(DA)$ if and only if $\lambda_N+\lambda_{\sigma(N)}=0$, for all $N\in\cal N$. From that the formulae follow at once.

\end{proof}

\begin{obser}
\label{remark}
In order facilitate the identification of the neat
equivalence classes, the following comments are helpful.
First observe that
the elementary relations of Definition 4.1 preserve the  circuit,
i.e., if $(\alpha ,\beta )\sim (\gamma ,\delta )$ then
$\overline{\alpha\beta}=\overline{\gamma\delta}$. Hence,
identification of neat equivalence classes can be done circuit by
circuit and, in particular, $r=\sum r_C$, with the sum indexed by
the set of circuits in $Q$ and $r_C$ being  the number of neat
equivalence classes $(\overline{\alpha ,\beta })$ such that
$\alpha\beta\in C$. On the other hand, each neat equivalence class
has a representative $(\alpha ,\beta )$ such that $\alpha a$ and
$a\alpha$ are zero paths, for all $a\in Q_1$. In particular, one
should only consider circuits $C$ containing a pair $(\alpha
,\beta )$ with the latter property. We illustrate this in the examples at the end of this section.
\end{obser}

\begin{coro}
Let $A=kQ/<Z>$ be a monomial algebra, where $k$ is a field of
characteristic different from $2$ and $Q$ is a finite connected
quiver. Then $\Alt (DA) =0$ if and only if the involution
considered in Lemma \ref{involution} is the identity.

If $k$ is a field of characteristic two, $\Alt (DA) =0$ if and only if
the set of neat cyclic pairs is empty.
\end{coro}

\begin{coro}
Let $A=kQ/<Q_2>$ be a connected $2$-nilpotent algebra . Then, unless
$Q$ is a loop,  we have
$$\dim \Alt (DA) =
\begin{cases}
\frac{\vert Q_1 \odot Q_1 \vert - \vert Q_1 \odot Q_0 \vert}{2} \ & \mbox{if $\car k \not =2$} \\
\frac{\vert Q_1 \odot Q_1 \vert + \vert Q_1 \odot Q_0 \vert}{2} \ & \mbox{if $\car k =2$.}
\end{cases}$$

In the case where $Q$ is a loop, $\Alt(DA)$ is zero  when $\car k \neq 2$,
and has dimension 2 when $\car k =2$.
\end{coro}

\begin{proof}
If the quiver $Q$ is a loop with vertex $u$ and
loop $a$, then $(\overline{a,u}) =
(\overline{u,a})$ and $(\overline{a,a})$ are the only neat
equivalence classes, both of which are symmetric. Then $r=s=2$ and the result
follows in this case. If $Q$ is not a loop then $Q_1 \odot Q_1$ is
the set of neat cyclic pairs, and the equivalence relation is
just the equality. In particular, $(\overline{a,b})$ is symmetric
if and only if $a=b$ is a loop. Then $r=\vert Q_1\odot Q_1 \vert $ and
$s=\vert Q_1\odot Q_0 \vert $ and we are done.
\end{proof}

Using \cite{cmrs}[Theorem 5.5], the combination of the formulae in
\cite{cs}[Theorem 1 and Proposition 2] and our Theorems
\ref{teorema1} (resp. \ref{char-p}) and \ref{teorema3} gives a
precise formula for the dimension of $HH^1(\Lambda )$, when
$\Lambda =TA$ is the trivial extension of the monomial algebra
$A$. We don't write down that formula in order to avoid excessive
technicalities. Our final result is a consequence of this (unwritten) formula:
we describe the
monomial algebras $A$ such that $\Lambda =TA$ has minimal
$HH^1(\Lambda )$, see Corollary \ref{minimal}. But, motivated by a question
pointed out by the referee, we present this result in as much generality as we know.
Recall that if $c= a_r^{\epsilon_r}...a_1^{\epsilon_1}$
is a walk in $Q$, with $a_i$ an arrow and
$\epsilon_i\in\{\pm 1\}$ for $i=1,...,r$, then the integer $\vert \sum_{1 \leq i \leq r}\epsilon_i \vert $ is called the {\it weight} of the walk $c$.
For the terminology used in the following proposition, see \cite{japs}.

\begin{prop}\label{homogeneo}
Let $Q$ be a finite oriented quiver, $A=kQ/I$ be a finite dimensional algebra
with $I$ homogeneous, and $\Lambda =TA$ be its trivial extension. The
following assertions are equivalent:

\begin{enumerate}
\item There is a $k$-linear isomorphism $Z(A)\cong HH^1(\Lambda )$.
\item $\dim HH^1(\Lambda )=1$.
\item $HH^1(A)=0$.
\item The abelianization $\bar{\Pi}$ of the fundamental group $\Pi =\pi_1(Q,I)$ is a finite group
of order coprime to $char(K)$, and $l(Q,I)\subseteq d(Q,I)$.
\end{enumerate}
If these conditions are satisfied, the quiver $Q$ has no oriented cycles and the weight of every closed walk in $Q$ is zero.
\end{prop}

\begin{proof}
$2)\Longrightarrow 1)$ is clear and $1)\Longrightarrow 3)$ follows
from \cite{cmrs}[Theorem 5.5].

$3)\Longrightarrow 4)$ It follows from  \cite{japs}[Corollary 4(2)].

$4)\Longrightarrow 2)$ If $c=a_r^{\epsilon_r}...a_1^{\epsilon_1}$
is a closed walk in $Q$, with $a_i$ an arrow and
$\epsilon_i\in\{\pm 1\}$ for $i=1,...,r$, then we put $\xi
(c)=\sum_{1\leq i\leq r}\epsilon_i$. In this way, we get a map
from the set of closed walks in $Q$ to the integers, which is
compatible with the relations of the group $\Pi=\pi_1(Q,I)$. Hence,
we get a map $\xi :\Pi\longrightarrow\mathbf{Z}$ which is clearly
a group homomorphism. Since $\bar{\Pi}$ is finite, $\Hom(\Pi
,\mathbf{Z})=0$ and in particular $\xi =0$. From that one gets
that $\xi (c)=0$ for every closed walk $c$. This is equivalent to
say that every closed walk in $Q$ has zero weight, and it implies that
$Q$ has no oriented cycles. Hence, one gets that $Z(A)=k$, $HH_1(A)=0$
(see \cite{c1}, \cite{ha} and \cite{c11}) and, since $\Alt_A(DA)$
is generated by cyclic pairs of paths, we also infer that
$\Alt_A(DA)=0$. To see that $HH^1(A)=0$ we have to prove condition (*) of \cite{japs}[Corollary 4(2)]. Indeed, every nilpotent derivation maps an arrow onto a linear combination of parallel paths of length greater than one. But the
existence of such parallel path would imply the existence of a
closed walk with nonzero weight. So the unique nilpotent derivation is the trivial one, and \cite{japs}[Corollary 4(2)] applies to give
that $HH^1(A)=0$.  Now assertion 2) follows from
\cite{cmrs}[Theorem 5.5].
\end{proof}

\begin{rem}
We do not know if conditions 2 and 3 might be equivalent in general for arbitrary finite dimensional algebras. The
referee has pointed out this question.
\end{rem}

The following result is immediate.
\begin{coro}
Let $A=kQ/I$ be a finite dimensional algebra with $I$ homogeneous such that $HH^1(A)=0$. Then $HH_1(A)=0$ and $\Alt_A(DA)=0$.
\end{coro}

Recall that if $Q$ is a finite quiver, then
$\bar{Q}$ is the graph obtained from $Q$ by forgetting the
orientation of the arrows.  In the particular case of monomial algebras we can replace condition 4) in Proposition \ref{homogeneo} by the following.

\begin{coro}\label {minimal}
Let $Q$ be a finite oriented quiver, $A=kQ/<Z>$ be
a monomial algebra and $\Lambda =TA$ be its trivial extension. The
following assertions are equivalent:

\begin{enumerate}
\item There is a $k$-linear isomorphism $Z(A)\cong HH^1(\Lambda )$.
\item $\dim HH^1(\Lambda )=1$.
\item $HH^1(A)=0$.
\item $Q$ is a quiver without  double arrows such that  $\bar{Q}$ is a
tree.
\end{enumerate}
\end{coro}

\begin{proof}
We only have to prove that for monomial algebras condition 4) is equivalent to the last condition in Proposition \ref{homogeneo}.  Indeed, if $A$ is monomial,
$\Pi =\pi_1(Q, Z)$ is the fundamental
group of the graph $\bar{Q}$, which is always free. Consequently,
$\bar{\Pi}$ is finite if and only if $\Pi =1$, which is equivalent to
condition 4).
\end{proof}

\begin{ejem}
\label{ejemplo}

Let us consider the algebra $A$ given as follows. Its quiver $Q$
is a crown with $Q_0=\mathbf{Z}_7$ and arrows $a_i:i\rightarrow
i+1$, for all $i\in Q_0$. The set of relations is
$Z=\{a_3a_2a_1a_0, a_4a_3a_2a_1, a_5a_4a_3a_2, a_2a_1a_0a_6a_5\}$.
In this case, the only relevant circuit is the canonical one
 and,  clearly,  all cyclic pairs in
$B\odot B$ are neat. To determine the neat equivalence classes,
observe that the only paths $\alpha$ such that $\alpha a$ and
$a\alpha$ are zero paths for all $a\in Q_1$, are exactly
$\{a_3a_2a_1, a_4a_3a_2, a_1a_0a_6a_5a_4a_3, a_2a_1a_0a_6\}$.  A
direct computation then shows that the elements of $\cal N$, i.e.,
the (neat) equivalence classes are
$\{(\overline{a_3a_2a_1,a_0a_6a_5a_4})=
(\overline{a_4a_3a_2,a_1a_0a_6a_5})=(\overline{a_2a_1a_0a_6,
a_5a_4a_3}), (\overline{a_1a_0a_6a_5a_4a_3, a_2})\}\ $\  and the
involution of Lemma \ref{involution} is just the transposition.
Hence $r=2$, $s=0$ and we have $\dim \Hom_{A-A}(DA,A)=2$ and $\dim
\Alt(DA)=1$, in any characteristic.

We also see that there are no non trivial strong circuits, that
the only efficient (not strong) circuit is the canonical one,
denoted by $C$ in the sequel, which is also a $p'$-circuit when $\car
k=p>0$. We have $w_C= \vert (Q_1\odot B)_C \vert - \vert (Q_0\odot
B)_C \vert =1-0=1$. Therefore $HH_1(A)=0$.

On the other hand, we have $Z(A)\cong k$ and $\dim
HH^1(A)=1$, the latter being easily deducible from
\cite{japs}[Corollary 2(3)] or \cite{cs}[Theorem 1]. As a result,
we conclude that $\dim HH^1(TA)=3$, independently of the
characteristic.

\end{ejem}

\medskip

\begin{ejem}

Let $A=kQ/<Z>$ be the algebra given by a crown $Q$ of length $n>1$,
that is, $Q_0=\mathbf{Z}_n$, $Q_1= \{ a_i:i\rightarrow i+1, \ \forall
i \in Q_0\}$, and take $Z= \{ a_0 \gamma^{m-1}\}$ where $\gamma = \beta a_0$ and $\beta = a_{n-1} \dots a_1$, for some $m>1$.

The set of non trivial circuits in $Q$ is
$\{ C_i = \overline{\gamma^i} \}_{i>0}$.  It is clear that $\{ C_i \}_{0<i<m}$ is the set of non trivial strong circuits, and $C_m$ is the unique efficient circuit.  Moreover, $C_m$ has multiplicity $m$, so it is an efficient $p'$--circuit in $Q$ if and only if $m$ is not divisible by $p$.  Now,
$w_{C_m}=\vert (Q_1 \odot B)_{C_m} \vert - \vert (Q_0 \odot B)_{C_m} \vert= 1-0=1$, so Theorems \ref{teorema1} and  \ref{char-p} imply that
$$\dim  HH_1(A) =
\begin{cases}
m  \ & \mbox{if $\car k = p>0 $ and $p$ divides $m$,} \\
m-1\ & \mbox{otherwise}.
\end{cases}$$

By Remark \ref{remark}, the neat equivalence classes are
represented by the cyclic pairs $( \gamma^{m-1} \beta, a_0
\gamma^i) $ for $i$ with $0 \leq i \leq m-2$, and they are all
symmetric and non equivalent. So $\dim \Alt(DA)=0$ if $\car k \neq
2$ and $\dim \Alt (DA) = m-1$ if $\car k =2$.

Using the terminology in \cite{cs}[Theorem 1] we have that  $\dim
HH^1(A) = \dim Z(A) - \vert Q_0 // B \vert + \vert Q_1 // B \vert
$, because all the elements in $Q_1 // B$ are glued and
admissible. By \cite{cs}[Proposition 2], we have $\dim Z(A)=m$.
Finally, one easily sees that $\vert Q_0 // B \vert =mn$, $ \vert
Q_1
// B \vert =mn-1$ and, hence, $\dim HH^1(A)= m- mn + mn -1 =m-1$.

We conclude that
$$\dim  HH^1(TA) =
\begin{cases}
4m-2\ & \mbox{if $\car k =2$ and $m$ is even,} \\
4m-3 \ & \mbox{if $\car k =2$ and $m$ is odd,} \\
3m-1 \ & \mbox{if $\car k =p > 2 $ and $m$ is divisible by $p$,} \\
3m-2 \ & \mbox{otherwise}.
\end{cases}$$

\end{ejem}

\vskip1cm \footnotesize \noindent C.C.:
\\D\'epartement de Math\'ematiques,
 Universit\'e de Montpellier
2,  \\F--34095 Montpellier cedex 5, France. \\{\tt Claude.Cibils@math.univ-montp2.fr}

\vskip3mm \noindent M.J.R:\\ Departamento de Matem\'atica, Universidad Nacional del Sur,\\Av. Alem 1253\\8000
Bah\'\i a Blanca, Argentina.\\ {\tt mredondo@criba.edu.ar}

\vskip3mm \noindent M.S.:\\ Departamento de Matem\'aticas, Universidad de Murcia,\\ Aptdo. 4021\\
30100 Espinardo, Murcia, Spain.\\ {\tt msaorinc@um.es}

\end{document}